\documentclass[reqno]{amsart}
\begin{document}
\def\N{{\mathbb N}}
\def\Z{{\mathbb Z}}
\def\R{{\rm I\!R}}
\def\GL{{\bf GL}}
\def\phi{\varphi}
\def\const{{\rm const}}

\def\L{{\mathcal L}}
\def\P{{\mathbb P}}
\def\Q{{\mathcal Q}}
\def\V{{\mathcal V}}
\def\Chi{K}
\def\Zeta{Z}
\def\1#1{{\widetilde{#1}}}

\def\C{{\mathchoice {\setbox0=\hbox{$\displaystyle\rm C$}\hbox{\hbox to0pt{\kern0.4\wd0\vrule height0.9\ht0\hss}\box0}}{\setbox0=\hbox{$\textstyle\rm C$}\hbox{\hbox to0pt{\kern0.4\wd0\vrule height0.9\ht0\hss}\box0}} {\setbox0=\hbox{$\scriptstyle\rm C$}\hbox{\hbox to0pt{\kern0.4\wd0\vrule height0.9\ht0\hss}\box0}} {\setbox0=\hbox{$\scriptscriptstyle\rm C$}\hbox{\hbox to0pt{\kern0.4\wd0\vrule height0.9\ht0\hss}\box0}}}}

\def\codim{{\rm codim}}
\def\tr{{\mathbf{Tr}}}
\def\M{{\mathcal M}}
\def\ov{\overline}
\def\m{\mapsto}
\def\rk{{\rm rank\,}}
\def\id{{\sf id}}
\def\Aut{{\sf Aut}}
\def\CR{{\rm CR}}
\def\gcd{{\rm gcd}}
\def\crd{\dim_{{\rm CR}}}
\def\crc{{\rm codim_{CR}}}
\def\eps{\varepsilon}
\def\gr{{\bf Gr}}
\def\gl{{\bf GL}}
\def\mod{{\rm mod\,}}
\emergencystretch9pt
\frenchspacing

\newtheorem{Th}{Theorem}[section]
\newtheorem{Def}{Definition}[section]
\newtheorem{Cor}{Corollary}[section]
\newtheorem{Satz}{Satz}[section]
\newtheorem{Prop}{Proposition}[section]
\newtheorem{Lemma}{Lemma}[section]
\newtheorem{Rem}{Remark}[section]
\newtheorem{Example}{Example}[section]

\title{Configurations of linear subspaces and rational invariants}
\author{Dmitri~Zaitsev}
\address{Mathematisches Institut, 
        Eberhard-Karls-Universit\"at T\"ubingen,
        72076 T\"ubingen, Germany,
        E-mail address: dmitri.zaitsev@uni-tuebingen.de}
\begin{abstract}
We construct a birational equivalence between certain quotients of $s$-tuples
of equidimensional linear subspaces of $\C^n$ and some quotients of 
products of square matrices modulo diagonal conjugations.
In particular, we prove the rationality of the quotient space of $s$-tuples of linear $2$-planes
in $\C^n$ modulo the diagonal $\gl_n(\C)$-action .
Furthermore, we compute generators of the field of the rational invariants explicitly.
\end{abstract}
\maketitle

\section{Introduction}

Let $\gr_{n,d}(\C)$ denote the Grassmannian of all $d$-dimensional linear subspaces in $\C^n$
and $\gl_n(\C)\times(\gr_{n,d}(\C))^s\to(\gr_{n,d}(\C))^s$ the canonical diagonal action.
I.~V.~Dolgachev posed the following question:

{\em Is the quotient $\gr_{n,2}(\C)^s/\gl_n(\C)$ (e.g. in the sense of Rosenlicht) always rational? }

Recall that a Rosenlicht quotient of an algebraic variety $X$ acted by an algebraic group $G$ is
an algebraic variety $V$ together with a rational map $X\to V$ whose generic fibers 
coincide with the $G$-orbits. 
Such quotients always exist and are unique up to birational isomorphisms (\cite{R}).
In the sequel all quotients will be assumed of this type.
An algebraic variety $Q$ is rational, if it is birationally equivalent to $\P^m$ with $m=\dim Q$.

We answer the above question positively by applying the rationality of the quotient $(\gl_2(\C))^2/\gl_2(\C)$,
where $\gl_2(\C)$ acts diagonally by conjugations (see~\cite{Pr67} and surveys~\cite{Do,Bru}):

\begin{Th}\label{2}
For all positive integers $n$ and $s$, the quotient $(\gr_{n,2}(\C))^s/\gl_n(\C)$ is rational.
Equivalently, the field of rational $\gl_n(\C)$-invariants on $(\gr_{n,2}(\C))^s$ is pure transcendental.
\end{Th}

The statement of Theorem~\ref{2} has been recently proved by Megyesi~\cite{Meg}
in the case $n=4$ and by Dolgachev~\cite{DoP} in the case of odd $n$.
Their proofs are independent of the present one.

More generally, we show the birational equivalence between $(\gr_{n,d}(\C))^s/\gl_n(\C)$
and certain quotients of matrix spaces.
Let $\gl_n(\C)\times (\gl_n(\C))^s \to (\gl_n(\C))^s$ be the action defined by
$(g,M_1,\ldots,M_s)\mapsto (gM_1g^{-1},\ldots,gM_sg^{-1})$. 
The first main result of the present paper consists of the following two statements:

\begin{Th}\label{main}
\begin{enumerate}
\item Let $s$ and $d$ be arbitrary positive integers and $n=rd$ for some integer $r>1$.
Then $(\gr_{n,d}(\C))^s/\gl_n(\C)$ is birationally equivalent to $(\gl_d(\C))^k/\gl_d(\C)$, where 
$$k=\begin{cases} (r-1)(s-r-1)& \text{ if } s>r+1\\
		  1& \text{ else }.\end{cases} $$
\item Let $s$ be arbitrary positive integer and $n=(2r+1)e$, $d=2e$ for some integers $r$ and $e$.
Then $(\gr_{n,d}(\C))^s/\gl_n(\C)$ is birationally equivalent to $(\gl_e(\C))^k/\gl_e(\C)$, where
$$k=\begin{cases} 2r-4+(4r-2)(s-r-2)& \text{ if } r>1,\, s>r+1 \\
		  2(s-4)& \text{ if } r=1,\,s> 4 \\ 
		  1& \text{ else } .\end{cases} $$
\end{enumerate}
\end{Th}

Theorem~\ref{2} follows from the first part of Theorem~\ref{main}
and from the rationality of the quotient of $(\gl_2(\C))^s$ by $\gl_2(\C)$ (see Procesi~\cite{Pr67})
in the case $n$ is even and from the second part of Theorem~\ref{main} in the case $n$ is odd.
More recently Formanek~\cite{Fo79,Fo80} proved the rationality of $(\gl_n(\C))^s/\gl_n(\C)$ for $n=3,4$.
Using this result together with Theorem~\ref{main} we obtain:

\begin{Th}
For every positive integer $s$, the following holds:
\begin{enumerate}
\item $(\gr_{n,3}(\C))^s/\gl_n(\C)$ is rational for $n = 0\, (\mod 3)$;
\item $(\gr_{n,4}(\C))^s/\gl_n(\C)$ is rational for $n = 0\, (\mod 2)$;
\item $(\gr_{n,5}(\C))^s/\gl_n(\C)$ is rational for $n = 3\, (\mod 6)$;
\item $(\gr_{n,8}(\C))^s/\gl_n(\C)$ is rational for $n = 4\, (\mod 8)$;
\end{enumerate}
\end{Th}

We refer the reader to~\cite{BS} for similar equivalences of stable rationalities
(an algebraic variety $V$ is stable rational, if $V\times\P^m$ is rational for some $m$).
We also refer to~\cite{GP71,GP72} for the classification of quadruples of linear subspaces 
of arbitrary dimensions and their invariants. 
In our situation, however, all rational invariants of the quadruples (i.e. the case $s=4$) are constant
unless $r=2$ in the first part of Theorem~\ref{main}.

Our method is based on constructing certain normal forms for our algebraic group actions.
We call an algebraic variety acted algebraically by an algebraic group $G$ a $G$-space.
A $G$-subspace is a $G$-invariant locally closed algebraic subvariety of $X$.
We use the standard notation $Gs:=\{gs:g\in G\}$ and 
$GS:=\{gs:g\in G, s\in S\}$, where $S\subset X$  is an arbitrary subset.

\begin{Def}
We say that $(S,H)$ is a {\em normal form} for $(X,G)$,
if $H\subset G$ is a subgroup and $S\subset X$ is a $H$-subspace
such that the following holds:
\begin{enumerate}
\item $GS$ is Zariski dense in $X$;
\item $Gs\cap S = Hs$ for all $s\in S$.
\end{enumerate}
\end{Def}

Clearly these conditions guarantee the birational equivalence of the quotients $X/G$ and $S/H$. 
In \S\ref{div} and ~\ref{most} we construct certain normal forms that are isomorphic
to the spaces of matrices as in Theorem~\ref{main}. Then in \S\ref{compute} we use
these normal forms for explicit computations of generators of the fields of rational invariants
in each case of Theorem~\ref{main}.

This method has also applications to biholomorphic automorphisms of nonsmooth bounded domains,
where the configurations of linear subspaces appear naturally as collections of tangent subspaces
to the so-called characteristic webs. We refer the reader to~\cite{Z98} for further details.

\section{Notation.}

For brevity we write $\gl_n$ and $\gr_{n,d}^s$
for $\gl_n(\C)$ and $(\gr_{n,d}(\C))^s$ respectively.
The actions on the products will be assumed diagonal unless otherwise specified.
For $E$ and $E'$ vector spaces, $\dim E \le\dim E'$, 
denote by $\gr_d(E)$ the Grassmannian of all $d$-dimensional subspaces of $E$,
by $\gl(E)$ the group of linear automorphisms of $E$ and
by $\gl(E_1,E_2)$ the space of all linear embeddings of $E_1$ into $E_2$.
A $G$-space $X$ is homogeneous (almost homogeneous),
if $G$ acts transitively on $X$ (on a Zariski dense subset of $X$).

\section{The case $d$ divides $n$.}\label{div}

In this section we study the diagonal $\gl_n$-action on $\gr_{n,d}^s$ with $n=rd$ for some integer $r$. 
Clearly the set of all $r$-tuples $(V_1,\ldots,V_r)\in\gr_{n,d}^r$
such that $\C^n=V_1\oplus\cdots\oplus V_r$ is Zariski open and $\gl_n$-homogeneous. 
This can be reformulated in terms of normal forms:

\begin{Lemma}\label{r}
Suppose that $r\ge 2$, $s\ge r$ and $(E_1,\ldots,E_r)\in\gr_{n,d}^r$
is such that $\C^n=E_1\oplus\cdots\oplus E_r$.
Define
$$S_1 = S_1(s):=\{(V_1,\ldots,V_s)\in \gr_{n,d}^s :
	(V_1,\ldots,V_r)=(E_1,\ldots,E_r) \} $$
and the group $H_1:= \gl(E_1)\times\cdots\times\gl(E_r)\subset\gl_n$.
Then $(S_1,H_1)$ is a normal form for $(\gr_{n,d}^s,\gl_n)$.
\end{Lemma}

Now fix the splitting $\C^n=E_1\oplus\cdots\oplus E_r$ as above.
Let $V\in\gr_{n,d}$ be such that its projection on each $E_i$,
$i=1,\ldots,r$, is bijective. In this case we say that $V$ is in {\em general position}
with respect to $(E_1,\ldots,E_r)$.
Clearly the subset of all $V$, which are in general position, is Zariski open in $\gr_{n,d}$.
Then for every $i=2,\ldots,r$,
the projection $V_i$ of $V$ on $E_1\times E_i$
is a graph of a linear isomorphism $\phi_i\colon E_1\to E_i$.
We claim that this correspondence between elements $V\in\gr_{n,d}$
and $(r-1)$-tuples of linear isomorphisms
$$(\phi_2,\ldots,\phi_r)\in \prod_{i=2}^r \gl(E_1,E_i)$$
is birational.

\begin{Lemma}\label{X1}
Under the assumptions of Lemma~\ref{r} define
$$
\Phi\colon \prod_{i=2}^r \gl(E_1,E_i) \to \gr_{n,d}, 
$$
$$
\Phi(\phi_2,\ldots,\phi_r) :=
	\{(z\oplus\phi_2(z)\oplus\cdots\oplus\phi_r(z)) : z\in E_1\}
$$
and the $H_1$-action on $\prod_{i=2}^r \gl(E_1,E_i)$ by
$$ ((g_1,\ldots,g_r),(\phi_2,\ldots,\phi_r)) \mapsto
	(g_2\circ\phi_2\circ g_1^{-1}, \ldots, g_r\circ\phi_r\circ g_1^{-1}).$$
Then $\Phi$ is birational and $H_1$-equivariant.
The image of $\Phi$ consists of all subspaces $V\in\gr_{n,d}$
which are in general position with respect to $(E_1,\ldots,E_r)$.
\end{Lemma}

\begin{proof}
Since $\Phi$ defines a system of standard coordinates
on $\gr_{n,d}$, it is birational.
The equivariance and the statement on the image are straightforward.
\end{proof}

Since $\prod_{i=2}^r \gl(E_1,E_i)$ is $H_1$-homogeneous, we obtain:

\begin{Cor}\label{trans}
Under the assumptions of Lemma~\ref{r}
let $S'_1(s)\subset S_1(s)$ be the Zariski open subset,
where $V_{r+1}$ is in general position with respect to $(E_1,\ldots, E_r)$.
Then $H_1$ acts transitively on $S'_1(r+1)$.
\end{Cor}

This proves statement~1 of Theorem~\ref{main} for $s\le r+1$.
Suppose now that $s > r+1$ and $E_{r+1}\in\gr_{n,d}$ is
in general position with respect to $(E_1,\ldots, E_r)$.

\begin{Lemma}\label{r+1}
Define $S_2=S_2(s):=\left\{(V_1,\ldots,V_s)\in S_1(s) : V_{r+1}=E_{r+1} \right\}$
and the group $H_2:=\gl(E_1)$ with the action on $S_2$ defined by the homomorphism
$$ H_2\to H_1, \quad g\mapsto
(g,\phi_2\circ g\circ \phi_2^{-1}, \ldots, \phi_r\circ g\circ \phi_r^{-1}),$$
where $(\phi_2,\ldots,\phi_r):=\Phi^{-1}(E_{r+1})$.
Then $(S_2, H_2)$ is a normal form for $(\gr_{n,d}^s$, $\gl_n)$.
\end{Lemma}

\begin{proof}
By Corollary~\ref{trans}, $H_1 S_2$ is Zariski open in $S_1$.
By straightforward calculations, $(S_2,H_2)$ is a normal form for $(S_1,H_1)$
and hence, by Lemma~\ref{r}, a normal form for $(\gr_{n,d}^s,\gl_n)$.
\end{proof}

Let $V\in\gr_{n,d}$ be one more subspace in general position
with respect to $(E_1,\ldots,E_r)$. Define $(\psi_2,\ldots,\psi_r):=\Phi^{-1}(V)$.
Clearly the map
\begin{equation}\label{oneto}
\prod_{i=2}^r \gl(E_1,E_i) \to \gl(E_1)^{r-1}, \quad
(\psi_2,\ldots,\psi_r)\mapsto (\phi_2^{-1}\circ\psi_2,\ldots,\phi_r^{-1}\circ\psi_r) 
\end{equation}
is one-to-one, where $(\phi_2,\ldots,\phi_r)=\Phi^{-1}(E_{r+1})$ is fixed.
Moreover, we obtain:

\begin{Lemma}\label{Hi}
Define
$\Chi\colon \gl(E_1)^{r-1} \to \gr_{n,d}$ by
$$\Chi(\chi_2,\ldots,\chi_r) :=
	\Phi(\phi_2\circ\chi_2,\ldots,\phi_r\circ\chi_r)$$
and the $H_2$-action on $X_2$ by
$$ (g,(\chi_2,\ldots,\chi_r)) \mapsto
	(g\circ\chi_2\circ g^{-1}, \ldots, g\circ\chi_r\circ g^{-1}).$$
Then $\Chi$ is birational and $H_2$-equivariant.
The image of $\Chi$ consists of all subspaces $V\in\gr_{n,d}$
which are in general position with respect to $(E_1,\ldots,E_r)$.
\end{Lemma}

\begin{proof}
The mapping $\Chi$ is birational as the composition of $\Phi$ and the map (\ref{oneto}).
The other statements are straightforward.
\end{proof}

\begin{Cor}\label{psi}
Under the assumptions of Lemma~\ref{r+1} define 
$$\Psi\colon S_2(s)\to (\gl(E_1)^{r-1})^{s-r-1},$$
$$\Psi(E_1,\ldots,E_{r+1},V_{r+2},\ldots,V_s):=
	(\Chi^{-1}(V_{r+2}),\ldots,\Chi^{-1}(V_s)). $$
Then $\Psi$ is birational and $H_2$-equivariant.
\end{Cor}

Finally we obtain together with Lemma~\ref{r+1}:

\begin{Cor}\label{ahom}
The space $\gr_{n,d}^s$, where $n=rd$,
is almost $\gl_n$-homogeneous if and only if $s\le r+1$.
If $s>r+1$, there is a normal form for $(\gr_{n,d}^s,\gl_n)$
which is isomorphic to $(\gl_d^{(r-1)(s-r-1)},\gl_d)$.
\end{Cor}

This implies the first part of Theorem~\ref{main}.

\section{The case $n=(2r+1)e$ and $d=2e$.}\label{most}

We start with an $r$-tuple $(V_1,\ldots,V_r)\in\gr_{n,d}^r$.
Clearly the subset of all $r$-tuples which form a direct sum
is Zariski open and $\gl_n$-homogeneous.
Fix an $r$-tuple $(E_1,\ldots,E_r)$ in this subset
and suppose that $s > r$. A straightforward calculation yields:

\begin{Lemma}\label{r'}
Define
$$S_1 = S_1(s):=\{(V_1,\ldots,V_s)\in \gr_{n,d}^s :
	(V_1,\ldots,V_r)=(E_1,\ldots,E_r) \} $$
and the group $H_1:= \{g\in\gl_n : g(E_i)=E_i \text{ for all } i=1,\ldots,r\}$
with the diagonal action on $S_1$.
Then $(S_1,H_1)$ is a normal form for $(\gr_{n,d}^s,\gl_n)$.
\end{Lemma}

Now we wish to parametrize the $d$-dimensional linear subspaces in $\C^n$
with respect to $(E_1,\ldots,E_r)$.
We say that $V$ is in {\em general position} with respect to
$(E_1,\ldots,E_r)$, if $\dim W=e$, where
$W:= V\cap (E_1\oplus\cdots\oplus E_r)$,
and the projection of $W$ on each $E_i$, $i=1,\ldots,r$,
is injective. Clearly the subset of all $V\in\gr_{n,d}$,
which are in general position, is Zariski open.

We first give another description of the subspaces
$W\in\gr_e(E_1\oplus\cdots\oplus E_r)$ which are in
general position with respect to $(E_1,\ldots, E_r)$.
For this let $A_i=A_i(V)\in\gr_e(E_i)$ be the projections of $W$
and let $\phi_i=\phi_i(V)\in \gl(A_1,E_i)$ be the linear isomorphisms
whose graphs are equal to the projections of $W$ on
$E_1\oplus E_i$, $i=2,\ldots,r$.
Denote by $X=X(E_1,\ldots,E_r)$ the space of all tuples $(A,\phi_2,\ldots,\phi_r)$, where $A\in\gr_e(E_1)$
and $\phi_i\in \gl(A,E_i)$, with the standard structure of a quasiprojective variety.

\begin{Lemma}\label{X1'}
Define
$$\Phi\colon X \to \gr_e(E_1\oplus\cdots\oplus E_r),$$
$$\Phi(A,\phi_2,\ldots,\phi_r) :=
	\{(z\oplus\phi_2(z)\oplus\ldots\oplus\phi_r(z)) : z\in A\}$$
and the $H_1$-action on $X$ by
$$ (g,(A,\phi_2,\ldots,\phi_r))\mapsto
	(g(A), g\circ \phi_2\circ g^{-1},\ldots, g\circ \phi_r\circ g^{-1}). $$
Then $\Phi$ is birational and $H_1$-equivariant.
The image of $\Phi$ consists of all subspaces $W\in X$
which are in general position with respect to $(E_1,\ldots,E_r)$.
\end{Lemma}

\begin{proof}
The proof is straightforward.
\end{proof}

\begin{Cor}\label{trans'}
Let $S'_1(s)\subset S_1(s)$ be the Zariski open subset
of all $s$-tuples $\V=(V_1,\ldots, V_s)$ such that
$V_{r+1}$ is in general position with respect to $(E_1,\ldots, E_r)$.
Then $H_1$ acts transitively on $S'_1(r+1)$.
\end{Cor}

\begin{proof}
It follows from the general position condition that
$$W(\V):=V_{r+1}\cap (E_1\oplus\cdots\oplus E_r) \in \gr_e(E_1\oplus\cdots\oplus E_r)$$
in the notation of Lemma~\ref{X1'}.
Let $\V,\V'\in S'_1(r+1)$ be arbitrary elements.
Since $X$ is $H_1$-homogeneous,
there exists $g_1\in H_1$ such that $g_1(W(\V))=W(\V')$.
Then there exists $g_2\in\gl_n$ with $g_2(V_{r+1})=V'_{r+1}$ and
$g_2|(E_1\oplus\cdots\oplus E_r) = \id$.
By the construction, $g_2\in H_1$ and the proof is finished.
\end{proof}

\begin{Cor}
For $s\le r+1$, $\gr_{n,d}^s$ is almost $\gl_n$-homogeneous.
\end{Cor}

Let $s > r+1$ and $E_{r+1}\in\gr_{n,d}$ be
in general position with respect to $(E_1,\ldots, E_r)$.
As a direct consequence of Corollary~\ref{trans'} we obtain:

\begin{Lemma}\label{r+1'}
Define $S_2:=\{(V_1,\ldots,V_s)\in S_1 : V_{r+1}=E_{r+1} \}$
and $H_2:=\{g\in\gl_n : g(E_i)=E_i \text{ for all } i=1,\ldots,r+1 \}$
with the diagonal action on $S_2$.
Then $(S_2,H_2)$ is a normal form for $(\gr_{n,d}^s,\gl_n)$.
\end{Lemma}

For the sequel we suppose that $(E_1,\ldots,E_{r+1})$ as above is fixed.
Define
$$(A,\phi_2,\ldots,\phi_r):=
	\Phi^{-1}(E_{r+1}\cap (E_1\oplus\cdots\oplus E_r)).$$
Let $V\in\gr_{n,d}$ be one more subspace
which is in general position with respect to $(E_1,\ldots,E_r)$.
Clearly $V$ is not uniquely determined by its $e$-dimensional
intersection $\Zeta_1(V):=V\cap(E_1\oplus\cdots\oplus E_r)$.
However, using $E_{r+1}$, we can consider another intersection
$\Zeta_2(V):=V\cap(E_1\oplus\cdots\oplus E_{r-1}\oplus E_{r+1})$.
Then the map
$$\Zeta:=(\Zeta_1,\Zeta_2)\colon \gr_{n,d}\to \gr_e(E_1\oplus\cdots\oplus E_r)\times \gr_e(E_1\oplus\cdots\oplus E_{r-1} \oplus E_{r+1})$$
is birational with the inverse
$\Zeta^{-1}\colon (W,W')\mapsto W+W'$.

Using the construction of $X$ and $\Phi$ for the tuple $(E_1,\ldots,E_{r-1},E_{r+1})$ instead of 
$(E_1,\ldots,E_r)$, we obtain $Y:= X(E_1,\ldots,E_{r-1},E_{r+1})$ and 
the $H_2$-equivariant birational map
$$\Psi\colon Y \to \gr_e(E_1\oplus\cdots\oplus E_{r-1} \oplus E_{r+1})$$
Combining $\Phi$, $\Psi$ and $\Zeta$ we obtain:

\begin{Cor}\label{Omega}
The composition
$\Omega:=\Zeta^{-1}\circ (\Phi,\Psi)\colon X\times Y\to \gr_{n,d}$
is birational and $H_2$-invariant.
\end{Cor}

In the following we consider the cases $r > 1$ and $r=1$ separately.

\subsection{ The case $r > 1$.}

For $V\in\gr_{n,d}$, set
\begin{equation}\label{B}
(B,\psi_2,\ldots,\psi_r):= \Phi^{-1}(V\cap (E_1\oplus\cdots\oplus E_r)) \in X,
\end{equation}
\begin{equation}\label{C}
(C,\chi_2,\ldots,\chi_{r-1},\chi_{r+1}):= \Psi^{-1}(V\cap (E_1\oplus\cdots\oplus E_{r-1}\oplus E_{r+1})) \in Y.
\end{equation}

In order to construct a smaller normal form,
fix $(A',\phi'_2,\ldots,\phi'_r)\in X$, $A''\in\gr_e(E_1)$
and $\chi'_{r+1}\in\gl(A'',E_{r+1})$ such that
$E_1=A\oplus A'$, $A''\cap A=A''\cap A'=\{0\}$,
$E_i=\phi_i(A)\oplus\phi'_i(A')$ for all $i=2,\ldots,r$,
and $\C^n=E_1\oplus\cdots\oplus E_r\oplus \chi'_{r+1}(A'')$.

\begin{Lemma}\label{S}
In the notation of (\ref{B},\ref{C}) define
$$S:=\{V\in\gr_{n,d} : B=A',\, C=A'',\,
	\psi_2=\phi'_2,\ldots,\psi_r=\phi'_r, \chi_{r+1}=\chi'_{r+1} \}.$$
Then $H_2S$ is Zariski open in $\gr_{n,d}$.
\end{Lemma}

\begin{proof}
For $V\in\gr_{n,d}$ generic, $E_1=A\oplus B$.
Hence there exists $g\in\gl(E_1)$ with $g|A=\id$ and $g(B)=A'$. 
Clearly $g$ extends to an isomorphism from the action by $H_2$.
Without loss of generality, $B=A'$.
Similarly there exists an isomorphism $g\in H_2$
such that $g|\phi_i(A)=\id$ and $g\circ\psi_i=\phi'_i$ on $A'$ for all $i=2,\ldots,s$.
Thus we may assume that $\psi_i=\phi'_i$.
By Lemma~\ref{trans}, we may also assume that $C=A''$.
Again, for $V$ generic, $\C^n=E_1\oplus\cdots\oplus E_r\oplus \chi_{r+1}(A'')$,
where $\chi_{r+1}(A'')\subset E_{r+1}$.
Therefore there exists an isomorphism $g\in H_2$ such that
$g|(E_1\oplus\cdots\oplus E_r)=\id$ and $g\circ\chi_{r+1}=\chi'_{r+1}$.
This proves the lemma.
\end{proof}

Let $\alpha\in\gl(A,A')$ be the isomorphism whose graph is $A''$
with respect to the splitting $E_1=A\oplus A'$.
Define $\beta\colon A\to A''$ by $\beta(z):=z\oplus\alpha(z)$.
For every $i=2,\ldots,r$, consider
$$\xi_i:=(\phi_i\times\phi'_i)^{-1}\circ\chi_i\circ\beta\in\gl(A,A\oplus A')$$
and its components
$[\xi_i]_1\in\gl(A)$, $[\xi_i]_2\in\gl(A,A')$.
Define the rational map
$$\Xi\colon S\to \gl(A)^{2(r-2)}, \quad \Xi\colon V\mapsto
	([\xi_2]_1,\alpha^{-1}\circ[\xi_2]_2,\ldots, [\xi_{r-1}]_1,\alpha^{-1}\circ[\xi_{r-1}]_2),$$
where $S$ is defined in Lemma~\ref{S}.

Let $g\in \gl(A)$ be arbitrary.
Using the isomorphisms $\alpha\colon A\to A'$,
$\phi_i\colon A\to \phi_i(A)$, $\phi'_i\colon A'\to \phi'_i(A')$, $i=2,\ldots,r$,
and $\chi'_{r+1}\circ\beta \colon A\to \chi'_{r+1}(A'')$
we can extend $g$ canonically to an isomorphism from the action by $H_2$.
This extension defines a canonical $\gl(A)$-action on $\C^n$ and therefore on $S$
(as defined in Lemma~\ref{S}).
We compare this action with the diagonal $\gl(A)$-action on $\gl(A)^{2(r-2)}$ by conjugations:

\begin{Lemma}\label{Xi}
The map $\Xi\colon \gr_{n,d}\to \gl(A)^{2(r-2)}$ is birational and $\gl(A)$-equivariant.
\end{Lemma}

\begin{proof}
The inverse of $\Xi$ is given by
$$\Xi^{-1}(\lambda_2,\ldots,\lambda_{r-1},\delta_2,\ldots,\delta_{r-1})=\Omega(x,y),$$
where $\Omega\colon X\times Y\to \gr_{n,d}$ is birational by Corollary~\ref{Omega} and
\begin{multline*}
x:=(A',\phi'_2,\ldots,\phi'_r), \quad
y:=(A'',(\phi_2\times\phi'_2)\circ (\lambda_2,\alpha\circ\delta_2)\circ\beta^{-1}, \ldots, \\
 (\phi_{r-1}\times\phi'_{r-1})\circ (\lambda_{r-1},\alpha\circ\delta_{r-1})\circ\beta^{-1},
	\chi'_{r+1} ).
\end{multline*}
The equivariance is straightforward.
\end{proof}

In the following corollary we use the space $S$ as defined in Lemma~\ref{S}.

\begin{Cor}
Define $S_3(s):=\{(V_1,\ldots,V_s)\in\gr_{n,d}^s : V_{r+2}\in S \}$
and consider the diagonal $\gl(A)$-action on $S_3$.
Then $(S_3,\gl(A))$ is a normal form for $(\gr_{n,d}^s,\gl_n)$.
\end{Cor}

\begin{Cor}
If $d=2e$ and $n=(2r+1)e$, 
the space $\gr_{n,d}^{r+2}$ is almost $\gl_n$-homogeneous
if and only if $r=2$.
\end{Cor}

Now let $V\in\gr_{n,d}$ be one more subspace.
Using the data fixed above we associate with $V$
a tuple of $4r-2$ linear automorphisms of $A$ as follows.
In the notation (\ref{B},\ref{C})
let $\tau_1\in\gl(A,A')$ and $\tau_2\in\gl(A,A')$
be the linear isomorphisms with the graphs $B$ and $C$ respectively.
As above, let $\alpha\in\gl(A,A')$ be the isomorphism whose graph is $A''$.
For $j=1,2$, set $\tau'_j(z):=z+\tau_j(z)$ ($\tau'_1\colon A\to B$, $\tau'_2\colon A\to C$) and
$\sigma_j:=\alpha^{-1}\circ\tau_j\in\gl(A)$, 
for $i=2,\ldots,r$, 
$$\zeta_i:=(\phi_i\times\phi'_i)^{-1}\circ\psi_i\circ\tau'_1 \in\gl(A,A\oplus A')$$
and for $i=2,\ldots,r-1$,
$$\theta_i:=(\phi_i\times\phi'_i)^{-1}\circ\chi_i\circ\tau'_2 \in\gl(A,A\oplus A').$$
We write $[\zeta_i]_1,[\theta_i]_1\in\gl(A)$ and
$[\zeta_i]_2,[\theta_i]_2\in\gl(A,A')$
for the corresponding components with respect to the splitting
$E_1=A\oplus A'$.
Then we obtain the linear automorphisms of $A$:
$x_i:=[\zeta_i]_1$, $y_i:=\alpha^{-1}\circ [\zeta_i]_2$
for $i=2,\ldots,r$ and
$z_i:=[\theta_i]_1$, $t_i:=\alpha^{-1}\circ [\theta_i]_2$
for $i=2,\ldots,r-1$.

In order to deal with the remainder term $\chi_{r+1}$ we
define $\phi_{r+1}\in\gl(A,\C^n)$ by
$$\phi_{r+1}(z):= z\oplus\phi_2(z)\oplus\cdots\oplus\phi_r(z).$$
It follows from the construction that $\phi_{r+1}(A)\in\gr_e(E_{r+1})$
and $E_{r+1}=\phi_{r+1}(A)\oplus\chi'_{r+1}(A'')$.
Using this splitting we define the remainder isomorphisms
$$\theta_r:=(\phi_{r+1}\times\chi'_{r+1})^{-1}\circ\chi_{r+1}\circ\tau'_2
	\in\gl(A,A\oplus A')$$
and $z_r:=[\theta_r]_1$, $t_r:=\alpha^{-1}\circ [\theta_r]_2$ $\in\gl(A)$.

\begin{Lemma}\label{Theta}
Let $\Theta\colon \gr_{n,d}\to\gl(A)^{4r-2}$ be given by
$$\Theta\colon V\mapsto (\sigma_1,\sigma_2,x_2,y_2,z_2,t_2,\ldots,x_r,y_r,z_r,t_r).$$
Then $\Theta$ is birational and $\gl(A)$-equivariant.
\end{Lemma}

\begin{proof}
In the above notation the inverse $\Theta^{-1}$ can be calculated as follows:
$$ \tau'_j:= \alpha\circ\sigma_j + \id,
 \quad B:= \tau'_1(A), \quad C:=\tau'_2(A),$$
$$ \psi_i=(\phi_i\times\phi'_i)\circ (x_i,\alpha\circ y_i) \circ (\tau'_1)^{-1}$$
for $i=2,\ldots,r$,
$$ \chi_i=(\phi_i\times\phi'_i)\circ (z_i,\alpha\circ t_i) \circ (\tau'_2)^{-1}$$
for $i=2,\ldots,r-2$,
$$ \chi_{r+1}=(\phi_{r+1}\times\chi'_{r+1})\circ (z_r,\alpha\circ t_r)
	\circ (\tau'_2)^{-1}$$
and finally
$$V=\Omega((B,\psi_2,\ldots,\psi_r),(C,\chi_2,\ldots,\chi_{r-1},\chi_{r+1})).$$
The equivariance is straightforward.
\end{proof}

\begin{Cor}\label{mostcor}
If $s\ge r+2$, the space $(S_3,\gl(A))$ is isomorphic to 
$$(\gl_e^{2r-4+(s-r-2)(4r-2)},\gl_e).$$
\end{Cor}

\begin{proof}
The required birational isomorphism is given by
\begin{equation}\label{connect}
(V_1,\ldots,V_s)\mapsto (\Xi(V_{r+2}),\Theta(V_{r+3}),\ldots,\Theta(V_s)),
\end{equation}
where $\Xi$ is birational by Lemma~\ref{Xi} and $\Theta$ is birational by Lemma~\ref{Theta}.
\end{proof}

This implies the second part of Theorem~\ref{main} in the case $r > 1$.

\subsection{ The case $r=1$.}\label{r1}

In this case we have $n=3e$.
Recall that we fixed $E_1,E_2\in\gr_{n,d}$ in general position
such that $\dim A=e$, where $A:=E_1\cap E_2$.
Choose another subspace $E_3\in\gr_{n,d}$
such that $\dim A_1=\dim A_2=e$,
where $A_j:=E_j\cap E_3$ for $j=1,2$.
Recall that we defined $H_2\subset \gl_n$ to be
the stabiliser of both $E_1$ and $E_2$.
Denote by $H_3\subset H_2$ the stabiliser of $E_3$.
Then $H_3=\gl(A)\times\gl(A_1)\times\gl(A_2)$
with respect to the splitting $\C^n=A\oplus A_1\oplus A_2$. 
The following is straightforward.

\begin{Lemma}
Suppose that $s\ge 3$ and define
$S_3=S_3(s):=\{(V_1,\ldots,V_s)\in S_2 : V_3=E_3 \}$.
Then $(S_3,H_3)$ is a normal form for $(\gr_{n,d}^s,\gl_n)$.
\end{Lemma}

Now we choose $B_j\in\gr_e(E_j)$
such that $B_j\cap A = B_j\cap A_j = \{0\}$ for $j=1,2$.
Then each $B_j$ can be seen as the graph of an isomorphism
$\phi_j\in\gl(A,A_j)$. On the other hand, for $V\in\gr_{n,d}$ generic,
the subspaces $C_j:=V\cap E_j$ are graphs of isomorphisms
$\psi_j\in\gl(A,A_j)$. Define 
$$g=(\id,\psi_1\circ\phi_1^{-1}, \psi_2\circ\phi_2^{-1})
	\in\gl(A\oplus A_1\oplus A_2).$$
Then $g\in H_3$ and $g(B_j)=C_j$ for $i=1,2$,
in particular, $V=B_1+B_2$. 
Together with Lemma~\ref{r+1} this proves the following:

\begin{Lemma}
Suppose that $s\ge 4$ and define
$$S_4=S_4(s):= \{(V_1,\ldots,V_s)\in S_3 : V_4\cap E_1=B_1, V_4\cap E_2=B_2 \}$$
and $H_4:=\gl(A)$. Define the $H_4$-action on $S_4$ via the homomorphism
$$H_4\to H_3, \quad g\mapsto 
	(g,\phi_1\circ g\circ \phi_1^{-1}, \phi_2\circ g\circ \phi_2^{-1}).$$
Then $(S_4,H_4)$ is a normal form for $(\gr_{n,d}^s,\gl_n)$.
\end{Lemma}

As a special case of Corollary~\ref{ahom}, we obtain:

\begin{Lemma}
The space $\gr_{3e,2e}^s$ is almost $\gl_{3e}$-homogeneous
if and only if $s\le 4$. If $s\ge 5$, there exists
a normal form which is isomorphic to $(\gl_e^{2(s-4)},\gl_e)$,
where $\gl_e$ acts diagonally by conjugations.
\end{Lemma}

This implies the second part of Theorem~\ref{main} in the case $r=1$.

\section{Computation of rational invariants.}\label{compute}

\subsection{The case $n=rd$.}\label{case1}
In order to compute the rational $\gl_{rd}$-invariants of $\gr_{rd,d}^s$ we represent
the elements of $\gr_{rd,d}^s$ by the equivalence classes of $rd\times ds$ matrices $M\in \C^{rd\times sd}$,
where the equivalence is taken under the right multiplication by $\gl_d^s$.
Then the diagonal $\gl_{rd}$-action on $\gr_{rd,d}^s$ corresponds 
to the left multiplication on $\C^{rd\times sd}$.
We start with $2d\times 2d$ matrices. Define
\begin{equation}\label{FF}
D\colon \gl_{2d}\to\gl_d, \quad
D
\begin{pmatrix}
A_{11} & A_{12} \\
A_{21} & A_{22}
\end{pmatrix} 
:=
A_{11}A_{21}^{-1}A_{22}A_{12}^{-1}.
\end{equation}

Then $D$ is a rational map that is invariant under the right multiplication by $\gl_d\times\gl_d$
and the left multiplication by $\{e\}\times\gl_d$, where $e\in\gl_d$ is the unit.
Moreover, $D$ is equivariant with respect to the left multiplication by $\gl_d\times\{e\}$.

More generally, let $M\in \C^{rd\times sd}$ be a rectangular matrix, where $r$ and $s$ are
arbitrary positive integers.
We split $M$ into $d\times d$ blocks $A_{ij}$, $i = 1,\ldots,r$, $j=1,\ldots,s$.
Then for every $i = 2,\ldots,r$, and $j=2,\ldots,s$, define
$$D_{ij}\colon \C^{rd\times sd} \to \gl_d, \quad
D_{ij}(M) := D 
\begin{pmatrix}
A_{11} & A_{1j} \\
A_{i1} & A_{ij}
\end{pmatrix}.
$$
Similarly to $D$, each $D_{ij}$ is a rational map that is invariant 
under the right multiplication by $\gl_d^s$ and under the left multiplication by $\{e\}\times\gl_d^{r-1}$.

Finally, we construct rational maps that are invariant under the left multiplication by the
larger group $\gl_{rd}$. For this, we assume $s > r + 1$, otherwise the left $\gl_{rd}$-multiplication
is almost homogeneous (see Corollary~\ref{ahom}) and, hence, all rational invariants are constant.
Then every matrix $M\in \C^{rd\times sd}$ from a Zariski open subset can be split into two blocks: 
$$M=(A B),\quad A\in \gl_{rd},\quad B\in \C^{rd\times (s-r)d}.$$
Define
$$\phi(M):= A^{-1}B\in \C^{rd\times (s-r)d}$$
and 
$$ G_{ij}(M):= D_{ij}(\phi(M)), \quad i=2,\ldots,r, \quad j=2,\ldots,s-r. $$

In particular, for $d=1$, $r=2$ and $s=4$, we obtain the classical double ratio:
$$
G_{22} 
\begin{pmatrix}
1   & 1   & 1   &   1 \\
z_1 & z_2 & z_3 & z_4
\end{pmatrix}
 = D 
\begin{pmatrix}
z_2-z_3 & z_2-z_4 \\
z_1-z_3 & z_1-z_4
\end{pmatrix} 
 = \frac{z_2-z_3}{z_1-z_3} : \frac{z_2-z_4}{z_1-z_4}.$$

By the invariance of $D_{ij}$, every $G_{ij}$ is invariant under the left multiplication by $\gl_{rd}$
and under the right multiplication by $\{e\}\times\gl_d^{s-1}$.
Moreover, due to the construction, $G_{ij}$ is equivariant with respect
to the subgroup $\gl_d\times\{e\} \subset \gl_d\times \gl_d^{s-1}$, 
where we take the right action on $\C^{rd\times sd}$
and the conjugation on the image space $\gl_d$. 
Therefore we obtain the following rational $\gl_{rd}\times \gl_d^s$-invariants:
\begin{equation}
I_{\alpha\beta} :=  \tr (G_{\alpha_1\beta_1}\cdot\ldots\cdot G_{\alpha_k\beta_k}) = 
\tr ((D_{\alpha_1\beta_1}\circ\phi)\cdot\ldots\cdot (D_{\alpha_k\beta_k}\circ\phi)),
\end{equation}
where $\alpha\in \{2,\ldots,r\}^k$ and $\beta\in \{2,\ldots,s-r\}^k$ are arbitrary multi-indices.

Using a result of Procesi~\cite{Pr76} 
we show that the invariant field is actually generated by these traces of monomials.

\begin{Th}
Let $d$, $r$, $s$ be arbitrary positive integers such that $s > r + 1$.
Then the field of rational $\gl_{rd}\times\gl_d^s$-invariants of $rd\times sd$ matrices is generated by 
the functions $I_{\alpha\beta}$, where $\alpha\in \{2,\ldots,r\}^k$ and $\beta\in \{2,\ldots,s-r\}^k$ and $k < 2^d$.
\end{Th}

\begin{proof}
Set $n:=rd$ as above. We write the spaces $E_1,\ldots,E_{r+1}, V_{r+2},\ldots, V_s\in\gr_{n,d}$ 
as in Lemma~\ref{r+1} in the form of an $rd\times sd$ matrix with $d\times d$ blocks as follow:

\begin{equation}\label{form0}
M := 
\begin{pmatrix}
          E & 0 & \cdots & 0 	  & E	   & E         & \cdots & E       \\
          0 & E & \cdots & 0 	  & E	   & V_{2,r+2} & \cdots & V_{2,s} \\
\vdots & \vdots & \ddots & \vdots & \vdots & \vdots    & \ddots & \vdots  \\
          0 & 0 & \cdots & E      & E	   & V_{r,r+2} & \cdots & V_{r,s} 
\end{pmatrix},
\end{equation}
where $E$ denotes the identity matrix of the size $d\times d$.

Let $\1{S_2} \subset \C^{rd\times sd}$ be subspace of all matrices of the form (\ref{form0}).
Clearly $E_1,\ldots,E_{r+1}$ fulfil the assumptions of Lemma~\ref{r+1}.
Then it follows from Lemma~\ref{r+1}, that 
$(\1{S_2},\gl_d)$ is a normal form for $(\C^{rd\times sd},\gl_{rd}\times\gl_d^s)$,
where $\gl_d$ acts on $\1{S_2}$ as the diagonal subgroup of $\gl_{rd}\times\gl_d^s$,
i.e. each $d\times d$ block is conjugated by the same matrix.

By the definition of a normal form, it is sufficient to prove that 
the invariant field of $(\1{S_2},\gl_d)$ is generated by the restrictions of $I_{\alpha\beta}$'s.
Let $M$ be given by (\ref{form0}). By the obvious calculations,
\begin{equation}\label{form}
\phi(M)= 
\begin{pmatrix}
          E	   & E         & \cdots & E       \\
          E	   & V_{2,r+2} & \cdots & V_{2,s} \\
          \vdots   & \vdots    & \ddots & \vdots  \\
          E	   & V_{r,r+2} & \cdots & V_{r,s} 
\end{pmatrix}
\end{equation}
and hence
\begin{equation}\label{mon}
I_{\alpha\beta}(M) = \tr (V_{\alpha_1,r+\beta_1}\cdot\ldots\cdot V_{\alpha_k,r+\beta_k}).
\end{equation}

By a theorem of Procesi~\cite{Pr76}, the polynomial invariants of $(r-1)(s-r-1)$-tuples 
$(V_{ij})$ with respect to the diagonal $\gl_d$-conjugations are generated by the monomials
of the form (\ref{mon}) with $k < 2^d$. Since all points of $(\1{S_2},\gl_d)$ are semi-stable
(see \cite{MFK93}), the categorical quotient $\1{S_2}//\gl_d$ exists and is given by these monomials.
Then the rational invariants on $(\1{S_2},\gl_d)$ are pullbacks of rational functions on $\1{S_2}//\gl_d$
and the proof is finished.
\end{proof}

\subsection{The case $n=3e$, $d=2e$.}
In this case the elements of $\gr_{n,d}^s$ are represented by the equivalence classes of matrices
$M\in \C^{3e\times 2se}$ with $3e\times 2e$ blocks $M_1,\ldots,M_s$.
As above we are looking for rational invariants with respect to left multiplications by $\gl_{3e}$
and right multiplications by $\gl_d^s$. We start with the case of $3$ blocks. 
Define the rational map
$$ \phi\colon \C^{3e\times 6e} \to \C^{3e\times 6e},$$
\begin{multline}\label{eq-phi}
\phi  
\begin{pmatrix}
A_1 & A_2 & A_3 \\
B_1 & B_2 & B_3
\end{pmatrix}
 := 
\begin{pmatrix}
E'          & E'          & E'       \\
B_1A_1^{-1} & B_2A_2^{-1} & B_3A_3^{-1}
\end{pmatrix}\\ =: 
\begin{pmatrix}
  E &   0 &   E &   0 &   E &   0        \\
  0 &   E &   0 &   E &   0 &   E        \\
c_1 & d_1 & c_2 & d_2 & c_3 & d_3
\end{pmatrix},
\end{multline}
where $A_1,A_2,A_3 \in \gl_d$, $B_1,B_2,B_3 \in \C^{e\times 2e}$
and $E'\in\gl_{2e}$, $E\in\gl_e$ denote the identity matrices.
We see the corresponding subspaces $E_1,E_2,E_3\in\gr_{n,d}$ as graphs
of linear maps given by the matrices $(c_1\,d_1)$, $(c_2\,d_2)$, $(c_3\,d_3)$.
Our first goal will be to compute the intersections $A:=E_1\cap E_2$ and $A_j:=E_j\cap E_3$, $j=1,2$.
For this, we set
$$ 
\begin{pmatrix}
x_1 & x_2 & x_3\\
y_1 & y_2 & y_3
\end{pmatrix}
:= \left(
\begin{pmatrix}
c_1 & d_1\\
c_2 & d_2
\end{pmatrix}^{-1}
,\,
\begin{pmatrix}
c_3 & d_3\\
c_1 & d_1
\end{pmatrix}^{-1}
,\,
\begin{pmatrix}
c_2 & d_2\\
c_3 & d_3
\end{pmatrix}^{-1}
\right)
\begin{pmatrix}
E\\
E
\end{pmatrix}
\in \C^{e\times 3e}.
$$
Then $x_1c_1+y_1d_1 = x_1c_2+y_1d_2 = E$ and therefore the $3e\times e$ matrix 
$$
\begin{pmatrix}
x_1\\
y_1\\
E
\end{pmatrix}
$$
represents the intersection $E_1\cap E_2$. Similarly the $3e\times e$ matrices
$$
\begin{pmatrix}
x_2\\
y_2\\
E
\end{pmatrix}
\quad\text{ and }\quad
\begin{pmatrix}
x_3\\
y_3\\
E
\end{pmatrix}
$$
represent the intersections $E_3\cap E_1$ and $E_2\cap E_3$ respectively.
Equivalently, the $3$-tuple $(E_1, E_2, E_3)$ can be represented by the matrix
\begin{equation}\label{1}
\begin{pmatrix}
x_1 & x_2 & x_3 & x_1 & x_2 & x_3\\
y_1 & y_2 & y_3 & y_1 & y_2 & y_3\\
E   & E   & E   & E   & E   & E
\end{pmatrix}
\end{equation}

Now take the general matrix $M=(M_1\ldots M_s)\in\C^{3e\times 2es}$.
Following the construction of \S\ref{r1} we bring $A$, $A_1$ and $A_2$,
i.e. the matrix (\ref{1}) to a normal form.
For this, consider the square matrix
$$
H(M) := 
\begin{pmatrix}
x_1 & x_2 & x_3\\
y_1 & y_2 & y_3\\
E   & E   & E
\end{pmatrix}
$$
and multiply $M$ by $H(M)^{-1}$:

\begin{equation}\label{long}
H(M)^{-1} M := \left(
\begin{array}{ccccccccccccccc}
E & 0 & 0 & E & 0 & 0 & C_{14} & D_{14} & \cdots & C_{1s} & D_{1s}\\
0 & E & 0 & 0 & E & 0 & C_{24} & D_{24} & \cdots & C_{2s} & D_{2s}\\
0 & 0 & E & 0 & 0 & E & C_{34} & D_{34} & \cdots & C_{3s} & D_{3s}\\
\end{array}
\right)
\end{equation}

Next we normalize the $3e\times 2e$ blocks:
$$
\begin{pmatrix}
C_{1i} & D_{1i}\\
C_{2i} & D_{2i}\\
C_{3i} & D_{3i}\\
\end{pmatrix} 
\begin{pmatrix}
C_{1i} & D_{1i}\\
C_{2i} & D_{2i}\\
\end{pmatrix}^{-1} =
\begin{pmatrix}
E            &                0 \\
0            &                E \\
\alpha_{2i-1}& \alpha_{2i}
\end{pmatrix} 
$$
and define
$$
\sigma_{2i-1}:=\alpha_{2i-1}\alpha_7^{-1}, \quad 
\sigma_{2i}:=\alpha_{2i}\alpha_8^{-1}, \quad i=5,\ldots,s.
$$

Comparing this construction with \S\ref{r1} we conclude that the matrices
$\sigma_j$, $j=9,\ldots,2s$, represent exactly the $2(s-4)$ matrices
in the isomorphic normal form $(\gl_e^{2(s-4)},\gl_e)$.
Using the theorem of Procesi~\cite{Pr76} and the arguments similar to \S\ref{case1} we obtain:

\begin{Th}
Let $e$ and $s$ be arbitrary positive integers such that $s > 4$.
Then the field of rational $\gl_{3e}\times\gl_{2e}^s$-invariants of $3e\times 2es$ matrices 
is generated by the functions 
$$J_\alpha:= \tr (\sigma_{\alpha_1}\cdot\ldots\cdot \sigma_{\alpha_k}) ,$$
where $\alpha\in \{9,\ldots,2s\}^k$ and $k < 2^e$.
\end{Th}

\subsection{The case $n=(2r+1)e$ and $d=2e$.}
Here an element $(V_1,\ldots,V_s)\in\gr_{n,d}^s$ is represented by a matrix 
$M\in \C^{(2r+1)e\times 2se}$ with $s$ blocks $M_1,\ldots,M_s$ of the size $(2r+1)e\times 2e$.
Again, the $\gl_n$-invariants on $\gr_{n,d}^s$ correspond 
to $\gl_n\times\gl_{2e}^s$-invariants on $\C^{(2r+1)e\times 2se}$.
Similarly to the previous paragraph, 
we start with a special case of a $(2r+1)e\times (2r+2)e$ matrix
and compute representatives of the intersections $E_{r+1}\cap (E_1\oplus\cdots\oplus E_r)$
and $E_{r+2}\cap (E_1\oplus\cdots\oplus E_r)$ and $E_{r+2}\cap (E_1\oplus\cdots\oplus E_{r-1}\oplus E_{r+1})$
as in \S\ref{most}. We start by normalizing the $n\times 2e$ blocks as in (\ref{eq-phi}):
\begin{equation}\label{eq-phi'}
\phi  
\begin{pmatrix}
A_1 & \ldots & A_s \\
B_1 & \ldots & B_s
\end{pmatrix}
:= 
\begin{pmatrix}
E'          & \ldots      & E'       \\
B_1A_1^{-1} & \ldots      & B_sA_s^{-1}
\end{pmatrix} =: 
\begin{pmatrix}
E'  & \ldots      & E'       \\
C_1 & \ldots      & C_s
\end{pmatrix}.
\end{equation}

In order to compute the first intersection $E_{r+1}\cap (E_1\oplus\cdots\oplus E_r)$ we consider
the corresponding system of linear equations:
$$
(E'\ldots  E')
\begin{pmatrix}
X_1 \\
\vdots\\
X_r
\end{pmatrix} 
= E' X_{r+1}, \quad
(C_1\,\ldots\, C_r)
\begin{pmatrix}
X_1 \\
\vdots\\
X_r
\end{pmatrix} 
= C_{r+1} X_{r+1},
$$
where each $X_i$ is a $2e\times e$ block.
Solving from the first equation  $X_{r+1} = X_1+\ldots + X_r$
and substituting this in the second we obtain:
$$
((C_1-C_{r+1})\,\ldots\, (C_r-C_{r+1}))
\begin{pmatrix}
X_1 \\
\vdots\\
X_r
\end{pmatrix} 
= 0.
$$
This is a system of $(2r+1)e\times e$ equations with $2re\times e$ variables.
Hence, for $C_1,\ldots,C_{r+1}$ in general position, it has a solution
which can be represented by rational $2e\times e$ block functions 
$$X_i = X_i(M), \quad i=1,\ldots,r.$$

Comparing this with the construction of \S\ref{most} we see that the blocks
$$
\begin{pmatrix}
E'  \\
C_i
\end{pmatrix}
(X_i),
\quad i=1,\ldots,r,
$$
represent the subspaces $A,\phi_2(A),\ldots,\phi_r(A)$.
For our normalization (Lemma~\ref{S}) we also need
the subspaces $A',\phi'_2(A'),\ldots,\phi'_r(A')$ and $\chi'_{r+1}$
which come from the other intersections. 
For them we can also solve the corresponding linear systems and
obtain rational $2e\times e$ block functions 
$$Y_i(M), \quad i=1,\ldots,r, \quad\text{ and }\quad Z_{r+1}(M),$$
such that the blocks
$$
\begin{pmatrix}
E'  \\
C_i
\end{pmatrix}
(Y_i(M)),
\quad i=1,\ldots,r,
\quad\text{ and }\quad
\begin{pmatrix}
E'  \\
C_{r+1}
\end{pmatrix}
(Z_{r+1}(M)),
$$
represent $A',\phi'_2(A'),\ldots,\phi'_r(A')$ and $\chi'_{r+1}$ respectively. 
As above we put these blocks together in the $n\times n$ matrix:
$$
H(M):= \left(
\begin{pmatrix}
E'  & \ldots      & E'       \\
C_1 & \ldots      & C_r
\end{pmatrix}
\begin{pmatrix}
X_1    & Y_1     & \cdots  & 0      & 0       \\
\vdots & \vdots  & \ddots  & \vdots & \vdots  \\
0      & 0       & \cdots  & X_r    & Y_r    
\end{pmatrix}, \,
\begin{pmatrix}
E'  \\
C_{r+1}
\end{pmatrix}
(Z_{r+1}(M))
\right).
$$
The $n\times e$ columns of $H(M)$ are exactly the representatives of 
$$A,\, A',\, \phi(A), \,\phi'(A'),\, \ldots,\, \phi_r(A),\, \phi'_r(A'),\, \chi'_{r+1}(A'')$$
in this order. If $M$ is in general position, $H(M)$ is invertible. As above, consider the matrix $H(M)^{-1}M =$
$$
\left(
\begin{array}{cccccccccccccccccccccccccc}
E      &0       &\cdots  &0      &0      &0      & a_1     & b_{1,r+2} & c_{1,r+2} &\cdots & b_{1,s} & c_{1,s}\\
0      &E       &\cdots  &0      &0      &0      & a_2     & b_{2,r+2} & c_{1,r+2} &\cdots & b_{1,s} & c_{1,s}\\
\vdots &\vdots  &\ddots  &\vdots &\vdots &\vdots & \vdots  & \vdots    &\vdots     &\ddots &\vdots   &\vdots  \\
0      &0       &\cdots  &E      &0      &0&a_{2r-1}&b_{2r-1,r+2}&c_{2r-1,r+2}&\cdots & b_{2r-1,s} & c_{2r-1,s}\\  
0      &0       &\cdots  &0      &E      &0&a_{2r}  &b_{2r,r+2}  &c_{2r,r+2} &\cdots & b_{2r,s} & c_{2r,s}\\   
0      &0       &\cdots  &0      &0      &E&a_{2r+1}&b_{2r+1,r+2}&c_{2r+1,r+2}&\cdots & b_{2r+1,s} & c_{2r+1,s}
\end{array}
\right).
$$

Furthermore, the property $E_{r+1}=\Phi(A,\phi_2,\ldots,\phi_r) + \chi'_{r+1}(A'')$ implies 
$$a_2=a_4=\cdots=a_{2r}=a_{2r+1}=0.$$
Using the property $W':=\Phi(A',\phi'_2,\ldots,\phi'_r)\in E_{r+2}$
we may assume that $W'$ is represented by the block
$$
\begin{pmatrix}
b_{1,r+2}\\
\vdots\\  
b_{2r+1,r+2} 
\end{pmatrix}
$$
and therefore 
$$b_{1,r+2}=b_{3,r+2}=\cdots=b_{2r+1,r+2}=0.$$

Furthermore, the matrix components of the map (\ref{connect}) can be calculated directly:

\begin{multline}\label{zzeta}
\Zeta(M) = \Bigg( 
D
\begin{pmatrix}
a_1 & c_{1,r+2}\\
a_3 & c_{3,r+2}
\end{pmatrix},\,
D
\begin{pmatrix}
b_{2,r+2} & c_{2,r+2}\\
b_{4,r+2} & c_{4,r+2}
\end{pmatrix},\,
\cdots, \\
D
\begin{pmatrix}
a_1 & c_{1,r+2}\\
a_{2r-3} & c_{2r-3,r+2}
\end{pmatrix},\,
D
\begin{pmatrix}
b_{2,r+2} & c_{2,r+2}\\
b_{2r-2,r+2} & c_{2r-2,r+2}
\end{pmatrix}
\Bigg)
\end{multline}
and

\begin{multline}\label{tteta}
\Theta_i(M) = \Bigg( 
D
\begin{pmatrix}
a_1 & b_{1,i}\\
a_3 & c_{3,i}
\end{pmatrix},\,
D
\begin{pmatrix}
b_{2,r+2} & c_{2,i}\\
b_{4,r+2} & c_{4,i}
\end{pmatrix},\,
\cdots,\,
D
\begin{pmatrix}
a_1 & b_{1,i}\\
a_{2r-1} & b_{2r-1,i}
\end{pmatrix},\\
D
\begin{pmatrix}
b_{2,r+2} & c_{2,r+2}\\
b_{2r,r+2} & c_{2r,r+2}
\end{pmatrix},\,
(c_{2r+1,r+2}^{-1}b_{2r+1,i}),\,(c_{2r+1,r+2}^{-1}c_{2r+1,i})
\Bigg),
\quad i=r+3,\ldots,s,
\end{multline}
where $D$ is the generalized double ratio defined by (\ref{FF}).

Comparing this with the proof of  Corollary~\ref{mostcor} and applying
the theorem of Procesi~\cite{Pr76} we obtain:

\begin{Th}
Let $e$, $r$ and $s$ be arbitrary positive integers such that $s \ge r+2$.
Then the field of rational $\gl_{(2r+1)e}\times\gl_{2e}^s$-invariants of $(2r+1)e\times 2es$ matrices 
is generated by the functions 
$$\tr (\sigma_1\cdot\ldots\cdot \sigma_k) ,$$
where each $\sigma_l$ is either a component of the map $\Zeta$ in {\rm(\ref{zzeta})} or
of one of the maps $\Theta_{r+3},\ldots,\Theta_s$ in {\rm(\ref{tteta})} and $k < 2^e$.
\end{Th}

%\bibliography{lit}
%\bibliographystyle{alpha}

%\end{document}

\end{document}